\documentclass[letterpaper, 10pt, conference]{ieeeconf}

\IEEEoverridecommandlockouts% Only needed for \thanks command
\usepackage[
	version=2017/03/05,
	abbr,
	algorithmic=noend,
]{clevethm}

\usepackage[
	version=2017/03/17,
	eqreset=false,
]{myPreamble}

\newcommand\range[2]{#1,\ldots,#2}
\renewcommandx\Fw[3][{1=\ell},{3=\gamma}]{#2-#3{\nabla}\!#1(#2)}%				forward step \Fw[f]{x}[\gamma]
\renewcommandx\FB[4][{1=\ell},{2=g},{4=\gamma}]{%								forward-backward step \FB[f][g]{x}[\gamma]
	\prox_{#4#2}\left(\Fw[#1]{#3}[#4]\right)
}%

\newcommand\algName{{\bf \texttt{PANOC}}}
\newcommand\tell{{\smash{\tilde\ell}}}
\newcommand\myref[1][]{\hyperref[alg:zerofpr2]{\algName}}

\pgfplotsset{compat=1.10}
\pgfplotsset{every axis/.append style={
                    label style={font=\scriptsize},
                    tick label style={font=\scriptsize}
                    }}

\newif\ifproofs
\proofstrue % either this or the following
\title{\LARGE\bf A Simple and Efficient Algorithm for Nonlinear Model Predictive Control}

\author{%
	Lorenzo~Stella,
	Andreas~Themelis,
	Pantelis~Sopasakis and
	Panagiotis~Patrinos%
	\thanks{%
		The authors are with the Department of Electrical Engineering (\textsc{esat-stadius}) and Optimization in Engineering Center (\textsc{optec}), KU Leuven, Kasteelpark  Arenberg 10, 3001 Leuven, Belgium.
		The first two authors are also affiliated with the IMT School for Advanced Studies Lucca, Piazza S. Francesco 17, 55100 Lucca, Italy.%
	}%
}

\begin{document}
\maketitle

\begin{abstract}
	We present \algName, a new algorithm for solving optimal control problems arising in nonlinear model predictive control (NMPC).
	A usual approach to this type of problems is sequential quadratic programming (SQP), which requires the solution of a quadratic program at every iteration and, consequently, inner iterative procedures.
	As a result, when the problem is ill-conditioned or the prediction horizon is large, each outer iteration becomes computationally very expensive.
	We propose a line-search algorithm that combines forward-backward iterations (FB) and Newton-type steps
	over the recently introduced forward-backward envelope (FBE), a continuous, real-valued, exact merit function for the original problem.
	The curvature information of Newton-type methods enables
%	neutralizes the known sensitivity to ill-conditioning and problem size of FB, yet maintaining the favorable global convergence properties.
%	In fact, the proposed algorithm achieves
	asymptotic superlinear rates under mild assumptions at the limit point, and the proposed algorithm is based on very simple operations: access to first-order information of the cost and dynamics and low-cost direct linear algebra.
	No inner iterative procedure nor Hessian evaluation is required, making our approach computationally simpler than SQP methods.
	The low-memory requirements and simple implementation make our method particularly suited for embedded NMPC applications.
\end{abstract}

% ~~~~~~~~~~~~~~~~~~~~~~~~~~~~~~~~~~~~~~~~~~~~~~~~~~~~~~~~~~~~~~~~~~~~~~~~~~~~ %
% ~~~~~~~~~~~~~~~~~~~~~~~~~~~~~~~~~~~~~~~~~~~~~~~~~~~~~~~~~~~~~~~~~~~~~~~~~~~~ %

\section{Introduction}
	 % \textcolor{red}{Lorenzo, Pantelis}

% Mention other approaches. SQP, IP, etc.
% Also algorithms for OC problems with simple constraints, Dunn, Bertsekas,\ldots
% Mention perhaps single vs multiple shooting.
%
% Paper by Cannon~\cite{cannon2004efficient}.
%
% Our contribution:\ldots
%
% Advantages
% \begin{itemize}
% \item $O(N)$ complexity per iteration, can handle large horizons
% \item compared to standard approaches $g$ can be nonconvex.
% \item no need for solving a QP
% \item leverages on advances in AD tools such as CasADi.
% \item No second-order information when quasi-Newton directions are employed
% \item When l-bfgs is used, minimal linear algebra operations are required. Just inner products.
% \item Future works. Further exploit the structure of the Newton system. We can further improve (future work, put into conclusions)
% \end{itemize}

Model predictive control (MPC) has become a popular strategy to implement feedback control loops for a variety of systems, due to its ability to take into account for constraints on inputs, states and outputs.
Its success is intimately tied to the availability of efficient, reliable algorithms for the solution of the underlying constrained optimization problem:
linear MPC requires solving a convex QP at every sampling step, for which the mature theory of convex optimization provides simple and robust methods with global convergence guarantees.

On the other hand, the vast majority of systems are nonlinear by nature, and nonlinear models often capture their dynamics much more accurately. For this reason nonlinear MPC (NMPC) is a well suited approach to design feedback controllers in many cases. At every sampling step, NMPC requires the solution of a general nonlinear program (NLP):
general approaches for NLP include sequential quadratic programming (SQP) and interior-point methods (IP)~\cite{nocedal2006numerical,bertsekas2016nonlinear}. Typically this NLP represents a discrete-time approximation of the continuous-time, and thus infinite-dimensional, constrained nonlinear optimal control problem, within a direct optimal control framework. Various ways exist for deriving a finite-dimensional NLP from a continuous-time optimal control problem, namely single shooting, multiple shooting and collocation methods, see e.g.~\cite{diehl2006fast,quirynen2017numerical}.
Although multiple shooting formulations (keeping the states as problem variables) are recently popular, single shooting formulations (implicitly eliminating the states) have traditionally been used to exploit the sequential structure in optimal control problems, see \cite{dunn1989efficient,wright1990solution} and~\cite[§2.6]{bertsekas2016nonlinear} for a textbook account.

%\lorenzo{Mention approaches to NMPC here and references; single- and multiple-shooting formulations.}

\subsection{Problems framework and motivation}

In this paper we deal with discrete-time, optimal control problems with nonlinear dynamics. This type of problems can be obtained, for example, by appropriately discretizing continuous-time problems.
Furthermore, we allow for nonsmooth (possibly nonconvex) penalties on the inputs: these can be (hard or soft) input constraints, or could be used for example to impose (group) sparsity on the input variables by using sparsity-inducing penalties. Note that problems with soft state constraints fit this framework by including an additional smooth penalty on the system state (\eg, the squared Euclidean distance from a constrained set), in the spirit of a generalized quadratic penalty method.

By eliminating the state variables and expressing the cost as a function of the inputs only (single-shooting formulation), the NMPC problems that we address can be reduced to the minimization of a smooth, nonconvex function \(f\) plus a nonsmooth (possibly nonconvex) penalty \(g\).
This is precisely the form of problems that can be solved by the proximal gradient method, also known as forward-backward splitting (FBS), see~\cite{attouch2013convergence}, a generalization of the projected gradient method.
FBS is a fixed-point iteration for solving a nonsmooth, nonlinear system of equations defining the stationary points of the cost function.
As such, its iterations are very simple and ideal for embedded applications.
However, the simplicity of FBS comes at the cost of slow convergence to stationary points.
In fact, like all first-order methods, the behaviour of FBS is greatly affected by the problem conditioning: in the case of NMPC, it is customary to have ill-conditioned problems due to the nonlinear dynamics and the horizon length.

\subsection{Contributions and related works}

We propose a new, simple method for solving NMPC problems.
The proposed algorithm is a line-search method for solving the fixed-point equations associated with FBS, using the so-called \emph{forward-backward envelope} (FBE) as merit function to determine the stepsize \cite{patrinos2013proximal,stella2016forward,themelis2016forward}.
We show that if the search directions are computed using quasi-Newton formulas, then the algorithm converges with superliner asymptotic rate to a stationary point.
Computing the directions and evaluating the FBE simply require the computation of the forward-backward mapping, therefore the proposed algorithm is based on the very same operations as FBS:
\begin{enumerate}
  \item evaluation of the gradient of the smooth cost, which is performed using automatic differentiation (AD);
  \item evaluation of the proximal mapping of the nonsmooth penalty, which usually has a very simple closed-form.
\end{enumerate}
In particular, no second-order information on the problem cost is required.
Toolboxes for AD that use code generation to evaluate gradients and Jacobians efficiently, such as CasADi \cite{andersson2012casadi}, are available.
Furthermore, limited-memory methods such as L-BFGS~\cite{liu1989limited} that only perform inner products can be used to determine line-search directions, making the algorithm completely matrix-free and well suited for embedded implementations and applications.

A similar approach was recently exploited to analyze and accelerate the convergence of another proximal splitting algorithm in nonconvex settings, namely the Douglas-Rachford splitting, and its dual counterpart ADMM \cite{themelis2017douglas}.

The paper is organized as follows: in \cref{sec:Formulation} we frame the family of problems which we target; in \cref{sec:Algorithm} we describe the proposed method and discuss its properties; \cref{sec:Simulations} shows numerical results with the proposed algorithm.

\section{Problem Formulation}
	\label{sec:Formulation}
	We consider the following finite-horizon
%optimal control
problem 
\begin{subequations}\label{eq:MPCprob}
\begin{align}
	\minimize{}
&
	{\textstyle\sum_{n=0}^{N-1}}\ell_n(x_n,u_n)+g_n(u_n)+\ell_N(x_N)
\\
	\stt{}
&
	\fillwidthof[r]{x_{n+1}}{x_0}=\bar x
\\
&
	x_{n+1}=f_n(x_n,u_n),\ n=\range{0}{N-1}
\end{align}
\end{subequations}
where \(\func{f_n}{\R^{n_x\times n_u}}{\R^{n_x}}\), \(n=\range{0}{N-1}\) are smooth mappings representing system dynamics \(\func{\ell_n}{\R^{n_x\times n_u}}{\R}\), \(n=\range{0}{N-1}\), and \(\func{\ell_N}{\R^{n_x}}{\R}\) are smooth functions representing stage and terminal costs respectively, and \(\func{g_n}{\R^{n_u}}{\Rinf}\), \(n=\range{0}{N-1}\), are possibly \emph{nonconvex, nonsmooth} and \emph{extended-real-valued} functions representing penalties on the inputs, \eg constraints.

We are interested in simple algorithms for~\eqref{eq:MPCprob}, \ie algorithms that do not involve a doubly iterative procedure, such as SQP methods.
One such algorithm is certainly forward-backward splitting (FBS), also known as the proximal gradient method.
Let \(\func{F}{\R^{Nn_u}}{\R^{(N+1)n_x}}\) be defined as
\[
	F(\range{u_0}{u_{N-1}})=(\range{F_0(u)}{F_N(u)}),
\]
where \(F_0 \equiv \bar x\), while
\[
	F_{n+1}(u)={}f_n(F_n(u),u_n),\quad n=\range{0}{N-1},
\]
and, denoting \(u=(\range{u_0}{u_{N-1}})\),
\begin{align*}
	\ell(u)={} & {\textstyle\sum_{n=0}^{N-1}} \ell_n(F_n(u),u_n)+\ell_N(F_N(u)),
\\
	g(u)={} & {\textstyle\sum_{n=0}^{N-1}} g_n(u_n).
\end{align*}
Then, problem \eqref{eq:MPCprob} can be expressed as
\begin{equation}\label{eq:Problem}
	\minimize_{u\in\R^{Nn_u}}\varphi(u)\equiv\ell(u)+g(u).
\end{equation}
The FBS scheme is based on simple iterations of the form
\begin{equation}\label{eq:FB}
	u^{k+1}
{}\in{}
	T_\gamma(u^k)
{}\coloneqq{}
	\FB{u^k},
\end{equation}
where $\gamma>0$ is a stepsize parameter. Here, $\prox_{\gamma g}$ is the (set-valued) proximal mapping of $g$:
\[
	\prox_{\gamma g}(u)
{}={}
	\argmin_{v\in\R^{Nn_u}}\set{g(v)+\tfrac{1}{2\gamma}\|v-u\|^2}.
\]
For instance, when $g$ is the indicator of a set the proximal mapping is the Euclidean projection onto the set. We assume that $g$ is simple enough so that the proximal mapping can be evaluated efficiently, and this is true in many examples.
%Notice that when $g$ is separable, as in the optimal control problem~\eqref{eq:MPCprob}, the update on the inputs can be executed in parallel via
%\[
%	u_n^{k+1}\in\prox_{\gamma g_n}(u_n^k-\gamma\nabla_{u_n}\ell(u)).
%\]

The gradient of \(\ell\) in \eqref{eq:FB} is efficiently calculated by backward automatic differentiation (also known as reverse mode AD, adjoint method, or backpropagation), see~\cite{dunn1989efficient}.
%as in \Cref{alg:BAD}~\cite{dunn1989efficient}.
%\begin{algorithm}
%	\caption{Backward AD}%
%	\label{alg:BAD}%
%	\input{TeX/Alg/BAD.tex}%
%\end{algorithm}

Iteration \eqref{eq:FB} is a direct extension of the usual gradient method for problems involving an additional nonsmooth term \(g\).
It is widely accepted in the optimization community that the gradient method can be very inefficient: in fact, for nonlinear optimal control problems where \(g=0\) (unconstrained optimal control problems) several more efficient algorithms have been proposed such as nonlinear conjugate gradient or Newton methods, see \cite[§2.6]{bertsekas2016nonlinear}.
However, when the additional nonsmooth term \(g\) is present, one is left with not many choices other than the proximal gradient method.
One option in the case of \(g=\indicator_C\), where \(C\) is a box, is to apply the two-metric projection method of Gafni \& Bertsekas \cite{gafni1984twometric}, the trust-region algorithm of \cite{lin1999newton}, or the limited-memory BFGS algorithm for bound constrained optimization in \cite{byrd1995limited}, or more generally, when \(C\) is a simple polyhedral set (one that is easy to project onto), the algorithms of \cite{calamai1987projected}.
When \(C\) has a more complicated structure, extensions of this class of methods become quite complex \cite{dunn1988projected}.
When \(g\) is a general nonsmooth (perhaps nonconvex) function (such as the sparsity inducing \(\ell_1\)-norm, the sum-of-\(\ell_2\)-norms to induce group sparsity, or the indicator of a nonconvex set such as a finite set of points) then the mentioned algorithms do not apply.

In the present paper we develop an algorithm that requires exactly the same computational oracle as FBS and thus fits embedded applications, but that exploits some curvature information about \eqref{eq:MPCprob} in order to converge much faster.

	\subsection{Handling state constraints}
		\label{sec:StateConstraints}
		The following more general problem allows to handle cases in which state variables are also subject to constraints:
%Many optimal control problems in practice include  constraints on the states. The situation is covered by the following more general formulation:
\begin{align*}
	\minimize{}
&
	{\textstyle\sum_{n=0}^{N-1}}
%	\!\Bigl[
		\ell_n(x_n,u_n)+g_n(u_n)+h_n(C_n(x_n,u_n))
%	\Bigr]
\\
&
	{}+\ell_N(x_N)+h_N(C_N(x_N))
\\
	\stt{}
&
	\fillwidthof[r]{x_{n+1}}{x_0}=\bar x
\\
&
	x_{n+1}=f_n(x_n,u_n),\quad n=\range{0}{N-1}
\end{align*}
where
%, additionally to problem \eqref{eq:MPCprob},
\(\func{h_n}{\R^{m_n}}{\Rinf}\), \(n=\range{0}{N}\) are proper, closed, convex functions with easily computable proximal mapping and \(\func{C_n}{\R^{n_x}\times\R^{n_u}}{\R^{m_n}}\), \(n\in\range{0}{N-1}\), and \(\func{C_N}{\R^{n_x}}{\R^{m_N}}\) are smooth mappings. For example, when $h_n$ are indicators of the nonpositive orthant then we are left with a classical state-constrained optimal control problem.

Next, consider
\(
	\func{G}{\R^{Nn_u}}{\R^{m_0}{\times}{\cdots}{\times}\R^{m_N}}
\)
defined as
\[
	G(\range{u_0}{u_{N-1}})=(\range{G_0(u)}{G_N(u)}),
\]
where \(G_n(u)=C_n(F_n(u),u_n)\) for \(n=\range{0}{N-1}\),
\(G_N(u)=C_N(F_N(u))\), and \(\func{h}{\R^{m_0}{\times}{\cdots}{\times}\R^{m_N}}{\Rinf}\) with
\[
	h(z)={\textstyle\sum_{n=0}^{N-1}}h_n(z_n)+h_N(z_n).
\]
The problem can now be expressed as
\[
	\minimize\ell(u)+g(u)+h(G(u)).
\]

A standard practice in MPC is to include state constraints in the cost function via penalties. The reason for doing so is to avoid ending up with an infeasible optimal control problem which can easily happen in practice due to disturbances and plant-model mismatch. The usual way of doing so is by relaxing state constraints using a quadratic penalty. Taking this approach one step further, we smoothen out \(h\)  by replacing it with its \emph{Moreau envelope} $h^{1/\mu}$, \emph{i.e.}, the value function of the parametric problem involved in the definition of the proximal mapping. Here $\mu$ acts as a penalty parameter: in the case of state constraints of the form $G(u)\in C$, one has
$h^{1/\mu}(G(u))= \tfrac{\mu}{2}\dist^2_C(G(u))$
and the larger the value of $\mu$, the larger is the penalty for violating the state constraints.

 It is well known that the Moreau envelope is smooth when $h$ is proper, closed, convex. In fact its gradient is given by
 $$\nabla h^{1/\mu}(z)=\mu(z-\prox_{h/\mu}(z)).$$
Since  \(G\) is also smooth (as the composition of smooth mappings), the following modified stage costs are smooth
\begin{align*}
	\tell_n(u)
{}={}&
	\ell_n(F_n(u),u_n)+h^{\nicefrac1\mu}_n(G_n(u)),\
	n=\range{0}{N{-}1}
\\
	\tell_N(u)
{}={}&
	\ell_N(F_N(u))+h^{\nicefrac1\mu}_N(G_N(u))
\end{align*}
and the same holds for the total cost, which we redefine as
\[
	\ell\gets{\textstyle\sum_{n=0}^N}\tilde\ell_n.
\]
Therefore soft-state-constrained problems have the same form \eqref{eq:Problem}.
\cref{alg:BADsoft} can be used to efficiently compute \(\nabla\ell\).
\begin{algorithm}[t]
	\caption{Backward AD (soft-constrained states)}
	\label{alg:BADsoft}
	\begin{algorithmic}[1]
\Inputs
	\(x_0\in\R^n\),~
	\(u=(\range{u_0}{u_{N-1}})\)%
%\Outputs
%	\(\ell(u)\),~
%	\(\nabla\ell(u)\)%
\STATE
	\(\ell(u)\gets 0\)%
\def\myVar#1{\fillwidthof[l]{x_{n+1}}{#1}}% % < [c], [r] or [l] for center, right or left alignment
\def\myEq{\fillwidthof[c]{{}\gets{}}{{}={}}}%
\FOR{\(n=\range{0}{N-1}\)}
	\STATE
		\(\myVar{s_n}\myEq\prox_{h_n/\mu}(C_n(x_n,u_n))\)%
	\STATE
		\(\myVar{q_n}\myEq\mu(C_n(x_n,u_n)-s_n)\)%
	\STATE
		\(\myVar{\ell(u)}\gets\ell(u)+\ell_n(x_n,u_n)+h(s_n)+\tfrac{1}{2\mu}\|q_n\|^2\)%
	\STATE
		\(\myVar{x_{n+1}}\myEq f_n(x_n,u_n)\)%
\ENDFOR
\STATE
	\(\myVar{s_N}\myEq\prox_{h_N/\mu}(C_N(x_N))\)%
\STATE
	\(\myVar{q_N}\myEq\mu(C_N(x_N)-s_N)\)%
\STATE
	\(\myVar{\ell(u)}\gets\ell(u)+\ell_N(x_N)+h_N(s_N)+\tfrac{1}{2\mu}\|q_N\|^2\)%
\STATE
	\(\myVar{p_N}\myEq\nabla_{x_N}\ell_N+\nabla_{x_N} C_N q_N\)%
%\STATE
%	\(\myVar{p_N}\myEq\nabla\ell_N(x_N)+\nabla C(x_N)q_N\)%
\FOR{\(n=\range{N-1}{0}\)}
\STATE
	\(p_n = \nabla_{x_n} f_n p_{n+1} + \nabla_{x_n}\ell_n + \nabla_{x_n}C_n q_n\)%
%\STATE
%	\(p_n{{}=}\hphantom{+{}}\nabla_{x_n} f_n(x_n,u_n)p_{n+1}+\nabla_{x_n}\ell_n(x_n,u_n)\)%
%\item[]%
%	\(\hphantom{p_n=}{+{}}\nabla_{x_n}C_n(x_n,u_n)q_n\)%
\STATE
	\(\nabla_{u_n}\ell(u) = \nabla_{u_n} f_n p_{n+1} + \nabla_{u_n} \ell_n + \nabla_{u_n}C_n q_n\)%
%\STATE
%	\(\nabla_{u_n}\ell(u){{}=}\hphantom{+{}}\nabla_{u_n} f_n(x_n,u_n)p_{n+1}+\nabla_{u_n} \ell_n(x_n,u_n)\)%
%\item[]%
%	\(\hphantom{\nabla_{u_n}\ell(u)=}{+{}}\nabla_{u_n}C_n(x_n,u_n)q_n\)%
\ENDFOR
\end{algorithmic}

\end{algorithm}
\begin{rem}\label{rem:muVector}
	We have considered the case where parameter \(\mu\) is a scalar for simplicity: \cref{alg:BADsoft} can immediately be adapted, in case \(h_n\) are separable, to the case where \(\mu\) is a vector of parameters, of dimension compatible with the separability structure of \(h_n\). Similarly, the penalty parameter \(\mu\) can be allowed to depend on \(n\).
\end{rem}

\section{Forward-Backward Newton Type Algorithm}
	\label{sec:Algorithm}
	First studied for convex problems, FB iterations \eqref{eq:FB} have been recently shown to converge for problems where both $\ell$ and $g$ are nonconvex \cite{attouch2013convergence}: if \(\ell\) has \(L_\ell\)-Lipschitz continuous gradient and \(g\) is lower-bounded, then for any \(\gamma\in(0,\nicefrac{1}{L_\ell})\) all accumulation points \(u^\star\) of sequences complying with \eqref{eq:FB} are \(\gamma\)-\emph{critical}, \ie they satisfy the condition
\[
	u^\star\in\FB{u^\star}.
\]
Moreover, if \(\ell+g\) is a Kurdyka-\L ojasiewicz function --- a mild property satisfied by all subanalytic functions, for instance --- then any bounded sequence \eqref{eq:FB} is globally convergent to a unique critical point.

Because of such favorable properties, and the fact that in many problems the proximal mapping is available in close form, FBS has been extensively employed and studied.
The downside of such simple algorithm is its slow tail convergence, being it \(Q\)-linear at best and with \(Q\)-factor typically close to one when the problem is ill-conditioned.
The employment of variable metrics, \eg coming from Newton-type schemes, can dramatically improve and robustify the convergence, at the cost of prohibitively complicating the proximal steps, which would require inner procedures possibly as hard as solving the original problem itself.

\subsection{Newton-type methods on generalized equations}
Instead of directly addressing the minimization problem, one could target the complementary problem of finding critical points by solving the inclusion (\emph{generalized equation})
\begin{equation}\label{eq:zerR}
	\text{find \(u^\star\) such that }
	0\in R_\gamma(u^\star) {}\coloneqq{} \tfrac 1\gamma\left[u-T_\gamma(u)\right],
\end{equation}
Here \(R_\gamma\)
%\[
%	R_\gamma(u)
%{}\coloneqq{}
%%	\tfrac 1\gamma\left(u-\FB u\right)
%	\tfrac 1\gamma\left[u-T_\gamma(u)\right]
%\]
is the (set-valued) \emph{fixed-point residual}.
Under very mild assumptions, \(R_\gamma\) is well-behaved close to critical points, and when close to a solution problem \eqref{eq:zerR} reduces to a classical equation, as opposed to generalized equation.
This motivates addressing the problem using Newton-type methods
\begin{equation}\label{eq:NT_FPR}
	u^{k+1}=u^k-H_kR_\gamma(u^k),
\end{equation}
where \(H_k\) are invertible operators that, ideally, capture curvature information of \(R_\gamma\) and enable superlinear or quadratic convergence when close enough to a solution.
In quasi-Newton schemes, \(H_k\) is a linear operator recursively updated so as to satisfy the (inverse) secant condition
\[
	u^{k+1}-u^k=H_{k+1}\bigl(R_\gamma(u^{k+1})-R_\gamma(u^k)\bigr),
\]
and under mild differentiability assumptions at a candidate limit point \(u^\star\), local superlinear convergence is achieved provided that the Dennis-Moré condition
\begin{equation}\label{eq:DM}
	\lim_{k\to\infty}
	\frac{
		\|R_\gamma(u^k)-JR_\gamma(u^\star)d^k\|
	}{
		\|d^k\|
	}
{}={}
	0
\end{equation}
is satisfied, where \(d^k=-H_kR_\gamma(u^k)\).

\subsection{Forward-backward envelope}
The drawback of iterations of the type \eqref{eq:NT_FPR} is that convergence can only be guaranteed provided that \(u^0\) is close enough to a solution.
%of any fast method for solving nonlinear equations encoding optimality conditions of an optimization problem,
%is
%globalization, namely
%the very fact of getting ``close enough'' to a solution, so as to converge.
In fact, without globalization strategies such type of methods are well known to even possibly diverge.
In~\cite{themelis2016forward} a globalization technique  is proposed, based on the \emph{forward-backward envelope}  (FBE)~\cite{patrinos2013proximal} (initially derived for convex problems~\cite{patrinos2014forward, stella2016forward}).
% In this paper we follow the \cref{sec:zerofpr} we adapt the proposed algorithm to
The FBE is an exact, continuous, real-valued penalty function for  \eqref{eq:Problem}, defined as
\begin{equation}\label{eq:FBE}
	\varphi_\gamma(u)
{}\coloneqq{}
	\ell(u)-\tfrac\gamma 2\|\nabla\ell(u)\|^2+g^\gamma(u-\gamma\nabla\ell(u)).
\end{equation}
%Under no additional assumptions, \(\varphi_\gamma\) enjoys the following favorable properties.
\begin{prop}[\cite{themelis2016forward}]\label{prop:FBE}
	For any \(\gamma>0\), \(\varphi_\gamma\) is a strictly continuous function satisfying
	\begin{enumerate}
		\item\label{prop:FBEleq}
			\(\varphi_\gamma\leq\varphi\);
		\item\label{prop:FBEgeq}
			\(\varphi(\bar u)\leq\varphi_\gamma(u)-\tfrac{1-\gamma L_\ell}{2\gamma}\|u-\bar u\|^2\) for any \(\bar u\in T_\gamma(u)\).
	\end{enumerate}
	In particular,
	\begin{enumerate}[resume]
		\item\label{prop:FBEeq}
			\(\varphi(u)=\varphi_\gamma(u)\) for any \(u\in\fix T_\gamma\);
		\item
			\(\inf\varphi=\inf\varphi_\gamma\) and \(\argmin\varphi=\argmin\varphi_\gamma\) for any \(\gamma<\nicefrac{1}{L_\ell}\).
	\end{enumerate}
	\ifproofs
	\begin{proof}
	From the definition \eqref{eq:FBE} and \cite[Ex. 10.32]{rockafellar2011variational}, it is apparent that \(\varphi_\gamma\) is strictly continuous.
	Moreover, applying the definition of the Moreau envelope \(g^\gamma\) and rearranging, we may rewrite \(\varphi_\gamma\) as
	\[
		\mathtight
		\varphi_\gamma(u)
	{}={}
		\inf_w\set{
			\ell(u)+\innprod{\nabla\ell(u)}{w-u}+g(w)+\tfrac{1}{2\gamma}\|w-u\|^2
		}
	\]
	and taking \(w=u\) in the expression above shows \ref{prop:FBEleq}.
	Observing that the set of minimizers \(w\) is exactly \(T_\gamma(u)\), for any \(\bar u\in T_\gamma(u)\) the quadratic upper bound on \(\ell\) \cite[Prop. A.24]{bertsekas2016nonlinear} implies
	\begin{align*}
		\varphi_\gamma(u)
	{}={} &
	% 	\smashoverbracket{
			\ell(u)
			{}+{}
			\innprod{\nabla\ell(u)}{\bar u-u}
	% 	}{
	% 		\geq\ell(\bar u)
	% 		{}-{}
	% 		\nicefrac{L_\ell}{2}
	% 		\|\bar u-u\|^2
	% 	}
		{}+{}
		g(\bar u)
		{}+{}
		\tfrac{1}{2\gamma}
		\|u-\bar u\|^2
	\\
	{}\geq{} &
		\varphi(\bar u)
		{}+{}
		\tfrac{1-\gamma L_\ell}{2\gamma}
		\|u-\bar u\|^2
	\end{align*}
	which is \ref{prop:FBEgeq}.
	The last two claims then follow.
	\end{proof}
	\else
	\begin{proof}
		See \cite{arXivVersion}.
	\end{proof}
	\fi
\end{prop}
By strict continuity, via Rademacher's theorem~\cite[Thm. 9.60]{rockafellar2011variational}, \(\nabla\ell\) and \(\varphi_\gamma\) are almost everywhere differentiable with
\[
	\nabla\varphi_\gamma(u)
{}={}
	Q_\gamma(u)R_\gamma(u),
\quad\text{where}\quad
	Q_\gamma(u)
{}\coloneqq{}
	I-\gamma\nabla^2\ell(u),
\]
see \cite{themelis2016forward}. Matrices \(Q_\gamma(u)\) are symmetric and defined for almost any \(u\); if \(\gamma<\nicefrac{1}{L_\ell}\), then \(Q_\gamma(u)\) is also positive definite wherever it exists.
If \(\ell\) is twice differentiable at a critical point \(u^\star\) and \(\prox_{\gamma g}\) is differentiable at \(u^\star-\gamma\nabla\ell(u^\star)\), then \(\varphi_\gamma\) is twice differentiable at \(u^\star\) with Hessian~\cite{themelis2016forward}
\begin{equation}\label{eq:FBEHess}
	\nabla^2\varphi_\gamma(u^\star)
{}={}
	Q_\gamma(u^\star)JR_\gamma(u^\star).
\end{equation}
A sufficient condition for \(\prox_{\gamma g}\) to comply with this requirement involves a mild property of prox-regularity and twice epi-differentiability, see \cite[§13]{rockafellar2011variational}.
\begin{thm}[Strong local minimality. \cite{themelis2016forward}]\label{thm:StrLocMin}
	Let \(\gamma<\nicefrac{1}{L_\ell}\) and suppose that \(\nabla\ell\) and \(\prox_{\gamma g}\) are differentiable at a critical point \(u^\star\) and at \(u^\star-\gamma\nabla\ell(u^\star)\), respectively.
	Then, \(u^\star\) is a strong local minimum for \(\varphi\) iff it is a strong local minimum for \(\varphi_\gamma\), in which case \(\nabla^2\varphi_\gamma(u^\star)\) is positive definite and \(JR_\gamma(u^\star)\) is invertible.
\end{thm}

\subsection{A superlinearly convergent algorithm based on FBS steps}\label{sec:zerofpr}
\begin{comment}
We will address \eqref{eq:Problem} or \eqref{eq:ProblemSoft} with \Cref{alg:zerofpr2}, a method based on the FBE and on the fast BFGS update directions.
Starting from a symmetric and positive definite matrix \(B\), a BFGS update returns another symmetric and positive definite matrix as follows
\begin{equation}\label{eq:BFGS}
% \mathtight
% 	B_{k+1}
% {}={}
% 	\begin{cases}[l@{~~}l]
% 		B_k + \frac{y_k\trans{y_k}}{\innprod{y_k}{s_k}} - \frac{B_ks_k\trans{(B_ks_k)}}{\innprod{s_k}{B_ks_k}}
% 		&
% 		\text{if~} \innprod{s_k}{y_k}>0
% 	\\
% 		B_k
% 		&
% 		\text{otherwise.}
% 	\end{cases}
%%
	B
{}\gets{}
	\begin{cases}[l@{~~}l]
		B + \frac{y\trans{y}}{\innprod{y}{s}} - \frac{Bs\trans{(Bs)}}{\innprod{s}{Bs}}
		&
		\text{if~} \innprod{s}{y}>0
	\\
		B
		&
		\text{otherwise.}
	\end{cases}
\end{equation}
\end{comment}
To the best of our knowledge, \cite{themelis2016forward} proposes the first algorithm with superlinear convergence guarantees that is entirely based on forward-backward iterations.
%Relying solely on such a simple black-box oracle as FBS is an extremely appealing property especially for embedded applications, where limited computational power might pose severe constraints to the operations which can be performed.
In this work, we propose \myref[] (Proximal Averaged Newton-type method for Optimal Control), a new linesearch method for problem \eqref{eq:Problem}, which is even simpler than the one of~\cite{themelis2016forward}, yet it maintains all the favorable convergence properties.
After a quick glance at the favorable properties of the FBE and its kinship with FBS, the methodology of the proposed scheme is elementary.
At each iteration, a forward-backward element \(\bar u^k\) is computed.
Then, a step is taken along a convex combination of the ``nominal'' FBS update direction \(-\gamma r^k\) and a candidate fast direction \(d^k\).
By appropriately \emph{averaging} between the two directions we can ensure sufficient decrease of the FBE, enabling global convergence.
When close to a solution, fast directions will take over and the iterations reduce to \(u^{k+1}=u^k+d^k\).
{\renewcommand\thealgorithm{\algName.}%
\begin{algorithm}[t]
	\caption{
%		Proximal Averaged Newton-type method for Optimal Control
	}%
	\label{alg:zerofpr2}%
	\begin{algorithmic}[1]
\Inputs
	\(\gamma\in(0,\nicefrac{1}{L_\ell})\),\
	\(\sigma\in(0,\gamma\tfrac{1-\gamma L_\ell}{2})\),\
	\(u_0\in\R^{Nn_u}\).%
% \item[]%
% 	\hspace*{0pt}\hphantom{\textbf{Outputs: }}%
% 	\(B_0\in\symm_{++}(\R^{Nn_u})\).%
%\Outputs
%	\(u^k\in\R^{Nn_u}\).%
\FOR{\(k=0,1,\ldots\)}
	\STATE
		Compute \(\nabla\ell(u^k)\) using
%		\Cref{alg:BAD} (or \ref{alg:BADsoft})%
		\cref{alg:BADsoft}
	\STATE\label{state:zerofpr2:FB}
		\(
			\bar u^k
		{}={}
			\prox_{\gamma g}
			\left(
				u^k-\gamma\nabla\!_{u}\ell(u^k)
			\right)
		\),~
		\(
			r^k
		{}={}
			\tfrac{u^k-\bar u^k}{\gamma}
		\)
 	\STATE\label{state:zerofpr2:d}
% 		{\bf if } \(k>0\),~
% 		update \(B_{k-1}\mapsto B_k\) as in \eqref{eq:BFGS} with \(s_{k-1}=u^k-u^{k-1}\) and \(y_{k-1}=r^k-r^{k-1}\)%
		Let \(d^k=-H_kr^k\) for some matrix \(H_k\in\R^{Nn_u\times Nn_u}\)%
	\STATE
		\(
			u^{k+1}
		{}={}
			u^k-(1-\tau_k)\gamma r^k+\tau_kd^k
		\),
		where \(\tau_k\) is the largest in \(\set{(\nicefrac12)^i}[i\in\N]\) such that
		\begin{equation}\label{eq:LS}
			\varphi_\gamma(u^{k+1})
		{}\leq{}
			\varphi_\gamma(u^k)-\sigma\|r^k\|^2
		\end{equation}
\ENDFOR{}%
\end{algorithmic}
\end{algorithm}
}%
The next results rigorously show these claims.
\begin{thm}[Global subsequential convergence]
Consider the iterates generated by \myref[].
Then, \(r^k\to 0\) square-summably, and the sequences \(\seq{u^k}\) and \(\seq{\bar u^k}\) have the same cluster points, all satisfying the necessary condition for local minimality \(u\in\FB[\ell]u\).
\ifproofs
\begin{proof}
	First, the algorithm is well defined, that is, the linesearch \eqref{eq:LS} always terminates in a finite number of backtrackings.
	In fact, since \(\sigma<\gamma\tfrac{1-\gamma L_\ell}{2}\) and \(u^{k+1}\to\bar u^k\) as \(\tau_k\searrow 0\), continuity of \(\varphi_\gamma\), \cref{prop:FBEgeq,prop:FBEleq} imply that for small enough \(\tau_k\) \eqref{eq:LS} holds.
% 	\[
% 		\varphi_\gamma(u^{k+1})
% 	{}\to{}
% 		\varphi_\gamma(\bar u^k)
% 	{}\leq{}
% 		\varphi_\gamma(u^k)
% 		{}-{}
% 		\gamma\tfrac{1-\gamma L_f}{2}\|r^k\|^2
% % 	{}<{}
% % 		\varphi_\gamma(u^k)
% % 		{}-{}
% % 		\sigma\|r^k\|^2
% 	\]
	Telescoping \eqref{eq:LS}, and since \(\inf\varphi_\gamma=\inf\varphi>-\infty\) we obtain that \(\sum_{k\in\N}\|r^k\|^2<\infty\), and in particular that \(r^k\to 0\).
	Suppose now that \(\seq{\bar u^k}[k\in K]\to u'\) for some \(u'\in\R^{Nn_u}\) and \(K\subseteq\N\).
	Then, since
	\(
		\|\bar u^k-u^k\|
	{}={}
		\gamma\|r^k\|
	{}\to{}
		0
	\),
	in particular \(\seq{u^k}[k\in K]\to u'\) as well.
	Similarly, the converse also holds, proving that \(\seq{u^k}\) and \(\seq{\bar u^k}\) have same cluster points.
	Moreover,
	\[
		\mathtight
		u^k
	\in
		\cball{\bar u^k}{\gamma\|r^k\|}
	\subseteq
		\FB[\ell][g\!\!\!]{u^k}
		+
		\cball{0}{\gamma\|r^k\|}
	\]
	and since
	\(
		\seq{
			u^k
			{}-{}
			\gamma\nabla\ell(u^k)
		}[k\in K]
	{}\to{}
		u'
		{}-{}
		\gamma\nabla\ell(u')
	\),
	from the outer semicontinuity of \(\prox_{\gamma g}\) \cite[Ex. 5.23(b)]{rockafellar2011variational} it follows that \(u'\in\FB{u'}\), concluding the proof.
\end{proof}
\else
\begin{proof}
	See \cite{arXivVersion}.
\end{proof}
\fi
\end{thm}
\begin{rem}[Lipschitz constant \(L_\ell\)]
In practice, no prior knowledge of the Lipschitz constant \(L_\ell\) is required for \myref[].
In fact, replacing \(L_\ell\) with an initial estimate \(L>0\), the following instruction can be added right after \cref{state:zerofpr2:FB}:

\noindent\hspace*{0pt}\vspace{-\baselineskip}%
\begin{algorithmic}[1]
\makeatletter%
	\renewcommand{\ALC@lno}{%
		\ALC@linenosize\arabic{ALC@line}bis\ALC@linenodelimiter%
	}%
\setcounter{ALC@line}{2}%
\makeatother%
\IF{
	\(\ell(\bar u^k) > \ell(u^k)-\gamma\tinnprod{\nabla\ell(u^k)}{r^k}+\tfrac L2\|\gamma r^k\|^2\)
}%
	\item[]%
		\(\gamma\gets\nicefrac \gamma2\),~
		\(L\gets 2L\),~
		\(\sigma\gets\nicefrac \sigma2\)~
		and go to \cref{state:zerofpr2:FB}.%
\ENDIF{}%
\end{algorithmic}
The above condition will fail to hold as soon as \(L\geq L_\ell\) \cite[Prop. A.24]{bertsekas2016nonlinear}, and consequently \(L\) is incremented only a finite number of times.
%Whenever the quadratic bound is violated with \(L\) in place of \(L_\ell\), the estimated Lipschitz constant \(L\) is increased and \(\gamma\) decreased accordingly.
%Since replacing \(L_\ell\) with any \(L\geq L_\ell\) still satisfies the quadratic Lipschitz upper bound \cite[Prop. A.24]{bertsekas2016nonlinear}, it follows that \(L\) is incremented only a finite number of times.
Therefore, there exists an iteration \(k_0\) starting from which \(\gamma\) and \(\sigma\) are constant, and all the results of the paper remain valid if such a strategy is implemented.

Moreover, since \(\bar u^k\in\dom g\) by construction, if \(g\) has bounded domain and the selected directions \(d^k\) are bounded (as it is the case for any ``reasonable'' implementation), it suffices that \(\nabla\ell\) is \emph{locally} Lipschitz-continuous (\ie, strictly continuous), and as such any \(\ell\in C^2\) would fit the requirement.
In fact, in such case all the sequences \(\seq{u^k}\) and \(\seq{\bar u^k}\) are contained in a compact enlargement \(\Omega\) of \(\dom\partial g\), and \(L_\ell\) can be then taken as \(\lip_\Omega(\nabla\ell)\), or adaptively retrieved in practice as indicated above.
This is the typical circumstance in (N)MPC where \(g\) encodes input constraints, which in realistic applications are bounded.
\end{rem}

Each evaluation of \(\varphi_\gamma\) in the left-hand side of the linesearch condition \eqref{eq:LS} requires one forward-backward step; \(\varphi_\gamma(u^k)\) on the right-hand side, instead, is available from the previous iteration.
In particular, in the best case of stepsize \(\tau_k=1\) being accepted, each iteration requires exactly one forward-backward step.
Under mild assumptions, this is the case when directions \(d^k\) satisfy the Dennis-Moré condition \eqref{eq:DM}, as shown in 
%\Cref{thm:Superlinear} later on.
the following result. This shows that the FBE does not prevent superlinear convergence of \myref[] when Newton-type directions are used: eventually, unit stepsize is accepted and \myref[] reduces to~\eqref{eq:NT_FPR}.
This is in stark contrast with the well known drawback of classical nonsmooth exact penalties (the so-called Maratos effect \cite[§5.3]{bertsekas2016nonlinear}).

%The Dennis-Moré condition is enjoyed by directions generated with quasi-Newton schemes under differentiability assumptions at the limit point.
%The rank-two updates of BFGS are the most widely used in smooth optimization, and thanks to the globalization guarantees of \myref[] they straightforwardly fit to our framework.
%Because of problem size, the limited-memory variant L-BFGS is mostly used, which does not require the computation or storage of full matrices \(H_k\) but simply keeps memory of a small number of last pairs \(s_k=u^{k+1}-u^k\) and \(y^k=r^{k+1}-r^k\), and retrieves \(d=-H_kr^k\) by simply performing scalar products.
%In \cref{sec:Simulations} we will show the efficiency of \myref[] with L-BFGS directions compared to plain FBS and state-of-the-art solvers.
%
%The next theorem shows that the FBE does not prevent superlinear convergence of \myref[] whenever Newton-type directions are used: eventually, unit stepsize is accepted and \myref[] reduces to~\eqref{eq:NT_FPR}.
%This is in stark contrast with the well known drawback of classical nonsmooth exact penalty functions (the so-called Maratos effect, see, \eg \cite[§5.3]{bertsekas2016nonlinear}).

\begin{thm}[Superlinear convergence]\label{thm:Superlinear}
Suppose that in \myref[] \(u^k\to u^\star\), for a strong local minimum \(u^\star\) of \(\varphi\) at which \(R_\gamma\) and \(\nabla\varphi_\gamma\) are strictly differentiable.
If \(\seq{H_k}\) satisfies the Dennis-Moré condition \eqref{eq:DM}, then \(\tau_k = 1\) is eventually always accepted and \(u^k\to u^\star\) at superlinear rate.%
\ifproofs
\begin{proof}
From \cref{thm:StrLocMin} we know that
\(
	G_\star
{}\coloneqq{}
	\nabla^2
	\varphi_\gamma(u^\star)
{}\succ{}
	0
\)
and that \(JR_\gamma(u^\star)\) is nonsingular.
Since \(\bar u^k\) and \(u^k\) converge to \(u^\star\), up to an index shifting we may assume that \(\seq{u^k}\) is contained in an open set in which \(\varphi_\gamma\) is differentiable and \(R_\gamma\) continuous.
Since \(r^k=R_\gamma(u^k)\to0\), from \eqref{eq:DM} it follows that \(d^k\to 0\).
Let \(u^{k+1}_0\coloneqq u^k+d^k\); by adding and substracting \(R_\gamma(u_0^{k+1})\) in the numerator of \eqref{eq:DM}, by strict differentiability of \(R_\gamma\) at \(u^\star\) we obtain
\begin{equation}\label{eq:Rzero}
	\lim_{k\to\infty}
	\frac{
		\|R_\gamma(u^{k+1}_0)\|
	}{
		\|d^k\|
	}
{}={}
	\lim_{k\to\infty}
	\frac{
		\|R_\gamma(u^{k+1}_0)\|
	}{
		\|u^{k+1}_0-u^k\|
	}
{}={}
	0.
\end{equation}
Since \(JR_\gamma(x^\star)\) is nonsingular and \(u^{k+1}_0\to u^\star\), there exists a constant \(\alpha>0\) such that
\(
	\|R_\gamma(u^{k+1}_0)\|
{}\geq{}
	\alpha
	\|u^{k+1}_0-u^\star\|
\)
for \(k\) large enough.
Combined with \eqref{eq:Rzero} we obtain
\[
	0
{}\leftarrow{}
	\frac{
		\|u^{k+1}_0-u^\star\|
	}{
		\|u^{k+1}_0-u^k\|
	}
{}\geq{}
	\frac{
		\|u^{k+1}_0-u^\star\|
	}{
		\|u^{k+1}_0-u^\star\|
		{}+{}
		\|u^k-u^\star\|
	}.
\]
Divinding numerator and denominator by \(\|u^k-u^\star\|\) yields
\begin{equation}\label{eq:lim2}
	\lim_{k\to\infty}{
		\frac{
			\|u^k+d^k-u^\star\|
		}{
			\|u^k-u^\star\|
		}
	}
{}={}
	\lim_{k\to\infty}{
		\frac{
			\|u^{k+1}_0-u^\star\|
		}{
			\|u^k-u^\star\|
		}
	}
{}={}
	0.
\end{equation}
Therefore,
\begin{align*}
	\varepsilon_k
{}\coloneqq{} &
	\frac{\varphi_\gamma(u^{k+1}_0) - \varphi_\gamma(u^\star)}{\varphi_\gamma(u^k) - \varphi_\gamma(u^\star)}
\\
{}={} &
	\frac{
		\frac12
		\innprod{G_\star (u^{k+1}_0-u^\star)}{u^{k+1}_0-u^\star}+o(\|u^{k+1}_0-u^\star\|^2)
	}{
		\frac12
		\innprod{G_\star(u^k-u^\star)}{u^k-u^\star}+o(\|u^k-u^\star\|^2)
	}
\\
{}\leq{} &
	\frac{
		\|G_\star\|
		\left(
			\frac{
				\|u^{k+1}_0-u^\star\|
			}{
				\|u^k-u^\star\|
			}
		\right)^2
		{}+{}
		\left(
			\frac{
				o(\|u^{k+1}_0-u^\star\|)
			}{
				\|u^k-u^\star\|
			}
		\right)^2
	}{
		\lambda_{\min}(G_\star)
		{}+{}
		\left(
			\frac{
				o(\|u^k-u^\star\|)
			}{
				\|u^k-u^\star\|
			}
		\right)^2
	}
	{}\to{}
	0
\end{align*}
as \(k\to\infty\).
Moreover, since \(\bar u^k\to u^\star\) and \(u^\star\) is a (strong) local minimum, eventually \(\varphi_\gamma(\bar u^k)\geq\varphi_\gamma(u^\star)\); combining with \cref{prop:FBEleq} we obtain
\[
	\varphi_\gamma(u^k)-\varphi_\gamma(u^\star)
{}\geq{}
	\varphi_\gamma(u^k)-\varphi_\gamma(\bar u^k)
{}\geq{}
	\gamma\tfrac{1-\gamma L_f}{2}\|r^k\|^2.
\]
Therefore,
\begin{align*}
	\varphi_\gamma(u^{k+1}_0) - \varphi_\gamma(u^k)
{}\leq{} &
	-(1-\varepsilon_k)\bigl(\varphi_\gamma(u^k) - \varphi_\gamma(u^\star)\bigr)
\\
{}\leq{} &
	-(1-\varepsilon_k)\gamma\tfrac{1-\gamma L_f}{2}\|r^k\|^2
\\
{}\leq{} &
	-\sigma\|r^k\|^2
\quad\text{for \(k\) large enough,}
\end{align*}
where the last inequality follows from the fact that \(\varepsilon_k\to 0\) and \(\sigma<\gamma\tfrac{1-\gamma L_f}{2}\), so that eventually \((1-\varepsilon_k)\gamma\tfrac{1-\gamma L_f}{2}\geq\sigma\).
Therefore, for large enough \(k\) the linesearch condition \eqref{eq:LS} holds with \(\tau_k=1\), and unitary step-size is always accepted.
In particular, the limit \eqref{eq:lim2} reads
\(\displaystyle
\lim_{k\to\infty}
	\tfrac{\|u^{k+1}-u^\star\|}{\|u^k-u^\star\|}
{}={}
	0
\),
proving \(\seq{u^k}\) to be superlinearly convergent.
\end{proof}
\else
\begin{proof}
	See \cite{arXivVersion}.
\end{proof}
\fi
\end{thm}

Strict differentiability of \(R_\gamma\) and \(\nabla\varphi\) does not require any smoothness condition on the nonsmooth function \(g\).
In fact, the required conditions hold as long as \(g\) is \emph{prox-regular} and has a generalized quadratic \emph{epigraphical Hessian} at the limit point; prox-regularity is a mild property enjoyed, for instance, by any convex function or a function whose effective domain is a discrete set, and similarly functions complying with the required generalized second-order properties are ubiquitous in optimization, see \cite{rockafellar2011variational,themelis2016forward} and references therein.
For example, \emph{partly smooth} functions are a comprehensive class of functions for which such properties hold; in fact, if the critical point \(u^\star\) satisfies the qualification \(-\nabla\ell(u^\star)\in\relint\partial g(u^\star)\) and \(g\) is prox-regular at \(u^\star\), then \(\prox_{\gamma g}\) is differentiable around \(\Fw{u^\star}\), see \cite{lewis2002active}.

The Dennis-Moré condition is enjoyed (under differentiability assumptions at the limit point) by directions generated with quasi-Newton schemes, the BFGS method being a prominent example.
%are the most widely used in smooth optimization, and thanks to the globalization guarantees of \myref[] they straightforwardly fit to our framework.
Because of the problems size, in \cref{sec:Simulations} we will show the efficiency of \myref[] with its limited-memory variant L-BFGS: this does not require storing the matrices \(H_k\), but instead keeps memory of a small number of pairs \(s_k=u^{k+1}-u^k\) and \(y^k=r^{k+1}-r^k\), and retrieves \(d=-H_kr^k\) by simply performing scalar products.
%In \cref{sec:Simulations} we will show the efficiency of \myref[] with L-BFGS directions compared to plain FBS and state-of-the-art solvers.

%
\section{Numerical Simulations}
	\label{sec:Simulations}
	
% Lorenzo, Pantelis
%
% \begin{itemize}
% \item chain of masses
% \item test RK
% \item soft state constraints
% \item comparison with multiple shooting formulations. Also SNOPT?
% \item if time permits, closed-loop simulations.
% \item sanity check for all the algorithms
% \end{itemize}

To test the efficacy of the proposed algorithm we consider a system composed of a sequence of masses connected by springs \cite{wirsching2006fast,vukov2013auto}.
The chain is composed by \(M\) masses: one end is connected to the origin, while a handle on the other end allows to control the chain.
Let us denote by \(p^i(t) \in \R^3\) the position of the \(i\)-th mass at time \(t\), for \(i =\range{1}{M+1}\), where \(p^{M+1}(t)\) is the position of the control handle.
The control action at each time instant is denoted as \(u(t) = \dot p^{M+1}(t) \in \R^3\), \ie we control the velocity of the handle.
Each body in the chain has mass \(m\), and the springs have constant \(D\) and rest length \(L\). By Hook's law we obtain the dynamics \cite{wirsching2006fast}:
\begin{align*}
	\ddot p^i {}={}& \tfrac{1}{m}(F_{i,i+1}-F_{i-1,i}) + a, \\
	F_{i,i+1} {}={}& D\left(1-\tfrac{L}{\|p^{i+1}-p^i\|}\right)(p^{i+1}-p^i).
\end{align*}
where \(a = (0, 0, -9.81)\) is the acceleration due to gravity.
Denoting \(v^i\) the velocity of mass \(i\), the state vector is
\[
	x(t) = (\range{p^1(t)}{p^{M+1}(t)}, \range{v^1(t)}{v^M(t)}).
\]
In the three-dimensional space there are \(n_x = 3(2M+1)\) state variables and \(n_u = 3\) input variables, and
\[
	\dot x = f_c(x, u) = (\range{v^1}{v^M}, u, \range{\ddot p^1}{\ddot p^M}).
\]

\subsection{Simulation scenario}
An equilibrium state of the system was computed with the control handle positioned at a given \(p_{\mathrm{end}}\in\R^3\). This was perturbed by applying a constant input \(u = (-1, 1, 1)\) for \(1\) second, to obtain the starting position of the chain. The goal is to drive the system back to the reference equilibrium state: this is achieved by solving, for \(T > 0\)
% (\textcolor{red}{Panos: Lorenzo, make the notation consistent with what we have before.})
\begin{equation}\label{eq:CostFunctional}\begin{aligned}
  \minimize_{x, u}\ & L_c(T) = {\textstyle\int_0^{T}} \ell_c(x(t),u(t))dt \\
	\stt\ & \dot x = f_c(x, u)
\end{aligned}\end{equation}
% where \(T\) is the prediction horizon, subject to the dynamics and input constraints
% \[\|u(t)\|_\infty \leq 1,\quad t\in[0,T].\]
% The cost function \(c\) in \eqref{eq:CostFunctional} is selected as
where
\begin{equation}
\label{eq:CostFunction}
% \nonumber
	\ell_c(x,u)
{}={}%&
	% \phantom{+{}}
	\beta\|p^{M+1}-p_{\mathrm{end}}\|_2^2 + \gamma{\textstyle\sum_{i=1}^M}\|v^i\|_2^2 + \delta\|u\|_2^2.
% \\
% &
	% +\smash{\sum_{i=1}^{M+1}}\tfrac{\mu_i}{2} \dist^2_S(x^i).
	% \vphantom{\sum_{i=1}}
\end{equation}

To discretize \eqref{eq:CostFunctional} we consider a sampling time \(t_{\mathrm{s}}\) such that \(T = Nt_{\mathrm{s}}\) and  piecewise constant input \(u\) accordingly: for \(n = 0,\ldots,N-1\),
\(u(t) = u_n\) for all \(t\in[nt_{\mathrm{s}}, (n+1)t_{\mathrm{s}})\).
Then
\( L_c(T) = {\textstyle\sum_{n=0}^{N-1}} {\textstyle\int_{nt_{\mathrm{s}}}^{(n+1)t_{\mathrm{s}}}} \ell_c(x(t),u_n)dt
\):
the problem is cast into the form \eqref{eq:MPCprob} by discretizing the integrals in the sum, and the system dynamics, by setting
\begin{subequations}\label{eq:sim:discreteCost}\begin{align}
		\ell_n(x_n,u_n)
	{}\approx{}&
		{\textstyle\int_{nt_{\mathrm{s}}}^{(n+1)t_{\mathrm{s}}}}\ell_c(x(t),u_n)dt, \\
		f_n(x_n,u_n)
	{}\approx{}&
		{\textstyle\int_{nt_{\mathrm{s}}}^{(n+1)t_{\mathrm{s}}}} f_c(x(t),u_n)dt,
\end{align}\end{subequations}
with the initial condition \(x(nt_{\mathrm{s}}) = x_n\), \(n=\range{0}{N-1}\).
Furthermore, we constrain the states and inputs by setting \(g_n\) and \(h_n\) as the indicator functions of the feasible sets as
\begin{align*}
		g_n(u)
	{}={}&
		\delta_{\|\cdot\|_\infty \leq 1}(u), \\
		h_n(C_n(x, u))
	{}={}&
		{\textstyle\sum_{i=1}^{M+1}}\delta_{\geq -0.1}(x^i_2).
\end{align*}
Since \(h_n\) is separable with respect to the different masses, we smoothen it by associating a parameter \(\mu_i\) to each component (see \cref{sec:StateConstraints} and \cref{rem:muVector} in particular):
\begin{equation}\label{eq:sim:SoftConstraints}
	h_n^{1/\mu}(C_n(x, u)) = {\textstyle\sum_{i=1}^{M+1}}\tfrac{\mu_i}{2}\left(\min\set{0, p^i_2 + 0.1}\right)^2.
\end{equation}

% The last summand in \eqref{eq:CostFunction} penalizes the \(M\) masses and the handle leaving a feasible region \(S\) of the space: this corresponds to smoothing \(\delta_S\), the indicator function of set \(S\), by replacing it with its Moreau envelope with parameter \(1/\mu_i\) as discussed in \cref{sec:StateConstraints}. In particular, we considered
% \[ S = \set{(x_1,x_2,x_3)\in\R^3}[x_2 \geq -0.1]. \]
%
% This problem can be approximated by a problem of the form \eqref{eq:MPCprob} by setting \(g = \delta_{\|\cdot\|_\infty \leq 1}\), the indicator function of the infinity-norm ball, and appropriately discretizing the integral in \eqref{eq:CostFunctional}.

In the simulations we have used \(T=4\) seconds and a sampling time \(t_{\mathrm s} = 0.1\) seconds, which gives a prediction horizon \(N=40\). Integrals \eqref{eq:sim:discreteCost} were approximated with a one-step 4th-order Runge-Kutta method.
We used CasADi \cite{andersson2012casadi} to implement the dynamics and cost function, and to efficiently evaluate their Jacobian and gradient.
The model parameters were set as \(M=5\), \(m=0.03\) (kg), \(D=0.1\) (N/m), and \(L = 0.033\) (m). In \eqref{eq:CostFunction} we set \(\beta = 1\), \(\gamma = 1\), and \(\delta = 0.01\). The coefficients for the soft state constraints \eqref{eq:sim:SoftConstraints} were set as
\(\mu_1 = \mu_2 = \mu_3 = 10^2\), \(\mu_4 = \mu_5 = \mu_6 = 10\).
% % The resulting configuration is shown in \Cref{fig:ChainPosition}.

% \begin{table}
%   \centering
%   \begin{tabular}{|r|l|ll|}
%     \hline
%     \multicolumn{1}{|c|}{Parameter} & \multicolumn{1}{c|}{Description} & \multicolumn{2}{c|}{Value} \\
%     \hline
%     \(M\) & number of masses & \(5\) & \\
%     \(m\) & mass & \(0.03\) & kg \\
%     \(D\) & spring constant & \(0.1\) & N/m \\
%     \(L\) & spring rest length & \(0.033\) & m \\
%     \hline
%   \end{tabular}
%   \caption{Parameters of the spring-mass system used in the numerical simulations.}\label{tbl:ModelParameters}
% \end{table}

% \begin{figure}
% 	\centering
% 	% \includetikz[width=0.4\textwidth]{chain_positions}
% 	\input{TeX/Tikz/chain_positions}
%   \caption{\textcolor{red}{Is this figure needed?} The initial position for the simulation, in the spring-mass system, was obtained perturbing a reference equilibrium state by appling a constant input for \(1\) second.}\label{fig:ChainPosition}
% \end{figure}

\subsection{Results}
\begin{figure}[tb]
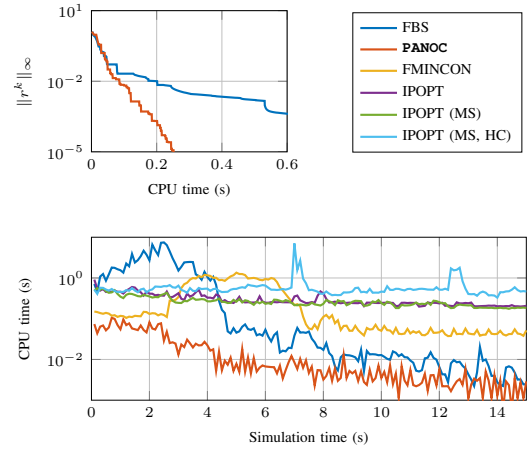

	\centering
	{{%
		\pgfkeys{/pgf/images/include external/.code={\includegraphics[width=.8\linewidth]{#width=.8\linewidth}}}%
		\tikzsetnextfilename{closed_loop_sim}%
		\input{./TeX/Tikz/closed_loop_sim.tex}%
	}}
	\caption{%
		(Top) Convergence of FBS and \myref[] in the first problem of the closed-loop simulation: the algorithms were executed here to medium/high accuracy for comparison purposes.
		(Bottom) CPU times of the solvers in the closed-loop simulation (``MS'': multiple shooting, ``HC'': hard constraints).%
	}%
	\label{fig:FPRConv}%
\end{figure}

\begin{figure}[tb]
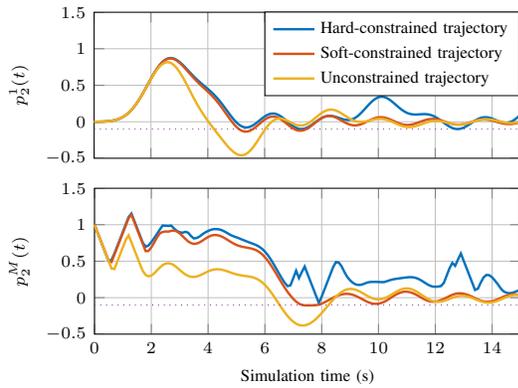

	\centering
	{{%
		\pgfkeys{/pgf/images/include external/.code={\includegraphics[width=.8\linewidth]{#width=.8\linewidth}}}%
		\tikzsetnextfilename{trajectory_1}%
		\input{./TeX/Tikz/trajectory_1.tex}%
	}}
	{{%
		\pgfkeys{/pgf/images/include external/.code={\includegraphics[width=.8\linewidth]{#width=.8\linewidth}}}%
		\tikzsetnextfilename{trajectory_M}%
		\input{./TeX/Tikz/trajectory_M.tex}%
	}}%
  \caption{Effect of the soft state constraint terms on the trajectory of the masses \(1\) and \(M\) in the closed loop simulation.
%  The unconstrained trajectory is obtained by setting \(\mu_i = 0\) for \(i=1,\ldots,M+1\) in \eqref{eq:sim:SoftConstraints}.
  }\label{fig:HandlePosition}
\end{figure}

We simulated the system for \(15\) seconds using different solvers.
In \myref{} we computed \(d^k\) in \cref{state:zerofpr2:d} using the L-BFGS method with memory \(10\) (see discussion in \cref{sec:zerofpr}).
Furthermore, we applied FBS, MATLAB's FMINCON (using an SQP algorithm), IPOPT (interior-point method) to both the single- and multiple-shooting formulations, and to the problem with hard state constraints.
We did not apply FMINCON to the multiple-shooting problem, as doing so performed considerably worse.
\cref{fig:FPRConv} shows the convergence of the fixed-point residual \(\|r^k\|_\infty\) for FBS and \myref[] for the first problem of the sequence: there we have solved the problem to medium/high accuracy for comparison purposes. In practice, we have noticed that good closed loop performance is obtained with more moderate accuracy: we ran closed-loop simulations terminating \myref[] and FBS as soon as \(\|r^k\|_\infty\leq 10^{-3}\). The other solvers were run with default options.
The CPU times during the simulation are shown in 
%\Cref{fig:ClosedLoopSim}.
\cref{fig:FPRConv}.
\myref[] outperforms the other considered methods in this example, and greatly accelerates over FBS:
this is particularly evident early in the simulation, when the system is far from equilibrium.
The effect of soft state constraints on the dynamics is shown in \cref{fig:HandlePosition}, where the trajectory of two masses during the simulation is compared to the hard-constrained and unconstrained cases.
Apparently, using soft state constraints improves considerably the solution time of the problem, without sacrificing closed loop performance.

%\subsection{Nonconvex input constraints}
%We have also applied FBS and \myref[] to the same problem with nonconvex input constraints. Specifically, with
%\[g_n(u) = \delta_{\|\cdot\|_2=1}(u), \]
%\ie, the indicator of the Euclidean unit sphere. The result of the closed-loop simulation is shown in \Cref{fig:ClosedLoopSimNoncvx}. Also here, like in the previous case, it is apparent that \myref[] solves the problem at each sampling step much faster than FBS.

% \begin{comment}
%\begin{figure}[tb]
%	\centering
%%	\includetikz[width=.76\linewidth]{closed_loop_sim}%
%	\input{TeX/Tikz/closed_loop_sim}
%	\caption{%
%		CPU times of the solvers applied to the spring-mass system, during \(15\) seconds of simulation.
%		``MS'' indicates that the multiple-shooting formulation was solved, where the system state is kept as decision variable. IPOPT was also applied to the problem with hard state constraints for reference.%
%	}%
%	\label{fig:ClosedLoopSim}%
%\end{figure}

\section{Conclusions}
	This paper presents \myref[], a new algorithm for solving nonlinear constrained optimal control problems typically arising in MPC. The algorithm is simple, exploits problem structure, does not require solution of a quadratic program at every iteration and yet can be shown to be superlinearly convergent under mild assumptions. Using L-BFGS directions in the algorithm was shown to perform favorably against state-of-the-art NLP solvers in a benchmark example.

There are several topics for future research:
\emph{(i)} semismooth Newton directions~\cite{patrinos2014forward}  that fully exploit the problem structure enabling quadratic convergence rates,
\emph{(ii)} more rigorous handling of state constraints by embedding the algorithm in a proximal augmented Lagrangian framework,
\emph{(iii)} a real-time iteration scheme where the algorithm is warm-started by exploiting sensitivity information for the fixed point residual and
\emph{(iv)} a code generation tool for embedded applications.

% ~~~~~~~~~~~~~~~~~~~~~~~~~~~~~~~~~~~~~~~~~~~~~~~~~~~~~~~~~~~~~~~~~~~~~~~~~~~~ %
% ~~~~~~~~~~~~~~~~~~~~~~~~~~~~~~~~~~~~~~~~~~~~~~~~~~~~~~~~~~~~~~~~~~~~~~~~~~~~ %

\bibliographystyle{IEEEtran}
\bibliography{TeX/Bibliography.bib}

\end{document}